\documentclass[11pt]{article}
\usepackage{amsmath}
\usepackage{amssymb}
\usepackage{graphicx}
\usepackage{tikz}
\usetikzlibrary{arrows.meta, positioning}

\numberwithin{equation}{section}

\newcommand{\FBOX}{\hspace*{\fill}$\rule{0.17cm}{0.17cm}$}
\newcommand{\BB}{\hspace*{\fill}$\rule{0.17cm}{0.17cm}\
\rule{0.17cm}{0.17cm}$}

\def\eref#1{(\ref{#1})}

\def\emeref#1{ {\em(\ref{#1})} } \def\emeref#1{ {\em(\ref{#1})}}

\newcommand{\Proof}{\noindent {\bf Proof.  }}

\newtheorem{THM}{THEOREM}[section] \newtheorem{LEMMA}[THM]{Lemma}
\newtheorem{CLAIM}[THM]{Claim} \newtheorem{COR}[THM]{Corollary}
\newtheorem{PROP}[THM]{Proposition}

\newtheorem{Thm}[THM]{Theorem} \newcommand{\Theorem}{\begin{Thm}}
\newcommand{\eTh}{\end{Thm}}

\newtheorem{CON}[THM]{Conjecture}
\newcommand{\Conjecture}{\begin{CON}} \newcommand{\eCon}{\end{CON}}

\newtheorem{CONS}[THM]{} \newcommand{\Conjectures}{\begin{CONS}}
\newcommand{\eCons}{\end{CONS}}

\newcommand{\eq}{\begin{equation}} \newcommand{\eeq}{\end{equation}}

\newcommand{\THEOREM}{\begin{THM}} \newcommand{\eT}{\end{THM}}
\newcommand{\Lemma}{\begin{LEMMA}} \newcommand{\eL}{\end{LEMMA}}
\newcommand{\Claim}{\begin{CLAIM}} \newcommand{\eCl}{\end{CLAIM}}
\newcommand{\Corollary}{\begin{COR}} \newcommand{\eCo}{\end{COR}}

\newcommand{\Proposition}{\begin{PROP}} \newcommand{\eP}{\end{PROP}}

\newtheorem{EXAM}{Example}[section]

\newcommand{\Example}{\begin{EXAM}} \newcommand{\eXa}{\end{EXAM}}

\newtheorem{EX}{Exercise}[section]

\newcommand{\Exercise}{\begin{EX}} \newcommand{\eEx}{\end{EX}}
\newtheorem{EXS}[EX]{} \newcommand{\Exercises}{\begin{EXS}}
\newcommand{\eExs}{\end{EXS}}

\newtheorem{PROBLEM}[EX]{Problem}
\newcommand{\Problem}{\begin{PROBLEM}}
\newcommand{\ePm}{\end{PROBLEM}}

\newtheorem{PROB}[EX]{} \newcommand{\Prob}{\begin{PROB}}
\newcommand{\eProb}{\end{PROB}}

\newtheorem{PROBLEMS}[EX]{} \newcommand{\Problems}{\begin{PROBLEMS}}
\newcommand{\ePms}{\end{PROBLEMS}}

\newtheorem{OPROBLEM}[EX]{Research problem}
\newcommand{\Open}{\begin{OPROBLEM}}
\newcommand{\eO}{\end{OPROBLEM}}

\newtheorem{REMARK}{Remark}[section]

\newcommand{\Remark}{\begin{REMARK}}
\newcommand{\eRe}{\end{REMARK}}

\newtheorem{COMMENT}{Comment}[section]
\newcommand{\Comment}{\begin{COMMENT}}
\newcommand{\eCom}{\end{COMMENT}}

\newtheorem{TODO}{Todo}[section] \newcommand{\Todo}{\begin{TODO}}
\newcommand{\eTo}{\end{TODO}}

\newtheorem{MEMO}{Memo}[section] \newcommand{\Memo}{\begin{MEMO}}
\newcommand{\eMe}{\end{MEMO}}

\newtheorem{ALGORITHM}[THM]{Algorithm}
\newcommand{\Algorithm}{\begin{ALGORITHM}}
\newcommand{\eAl}{\end{ALGORITHM}}

\begin{document}

\title{A simplified min-max formula \break for the inverse
arborescence problem}

\author{Andr\'as Frank \thanks{ELKH-ELTE Egerv\'ary Research Group,
Department of Operations Research, E\"otv\"os Lor\'and University
Budapest, P\'azm\'any P. s. 1/c, Budapest H-1117.  E-mail:  {\tt
andras.frank\char'100 ttk.elte.hu } } \ \ \ and \ \ \ {Hanna Szabrina
Horv\'ath \thanks{Department of Operations Research, E\"otv\"os
Lor\'and University Budapest, P\'azm\'any P. s. 1/c, Budapest H-1117.
E-mail:  {\tt horvath.hanna.szabrina\char'100 gmail.com }}} }

\maketitle

\begin{abstract} A simple min-max theorem is formulated and proved for
the smallest modification (measured in $l_1$-norm) of an input cost
function $w_0$ that makes a target arborescence $F_0$ of a digraph a
cheapest arborescence.  The constructive proof gives rise to a
polynomial time algorithm for computing both a minimizer cost function
on the primal side and a maximizer dual object in the min-max formula.
\end{abstract}


\section{Introduction}

A digraph $D=(V,A)$ with a specified root-node $r_0\in V$ is called
{\bf root-connected} if there is a dipath from $r_0$ to every other node.
This is equivalent to requiring that $D$ has a spanning arborescence
of root $r_0$ which will be abbreviated as an arborescence or an
$r_0$-arborescence.  More generally, $D$ is called {\bf rooted
$k$-arc-connected} if there are $k$ arc-disjoint dipaths from $r_0$ to
every other node.  By the directed arc-version of Menger's theorem,
$D$ is rooted $k$-arc-connected if and only if $\varrho _D(X)\geq k$
for every non-empty subset $X\subseteq V-r_0$, where $\varrho _D(X)$
denotes the number of arcs entering $X$.

Let $X\oplus Y :  = (X-Y) \cup (Y-X)$ denote the symmetric difference
of sets $X$ and $Y$.  For a function $w:S\rightarrow {\bf R} $, the
set function $\widetilde w$ is defined by $\widetilde w(X):= \sum
[w(s): \ s\in X]$ where $X\subseteq S$.

\subsection{ Classical results on arborescences }

Let $c:A\rightarrow {\bf R}_+$ be a non-negative cost function on the
arc-set of $D$.  We say that a function $y:  {\cal F}\rightarrow {\bf
R}$ defined on a set-system ${\cal F}\subseteq 2\sp V$ is {\bf
$c$-feasible} if $y$ is non-negative and \eq \hbox{ $\sum [y(Z): \ Z\in
{\cal F}, \ Z$ \ is entered by \ $a]\leq c(a)$ \ for every arc $a\in
A$.  }\ \label{(feas-def)} \eeq When ${\cal F}:=\{X: \ \emptyset \not =
X \subseteq V-r_0\}$, a $c$-feasible function $y$ will be referred to
as a {\bf dual solution} to the cheapest $r_0$-arborescence problem.
We call an arc $a$ of $D$ {\bf $c$-tight} (or just {\bf tight}) (with
respect to $y$) if $\sum [y(Z): \ Z\in {\cal F}, \ Z$ \ is entered by \
$a] = c(a)$.  Bock \cite{Bock} and Fulkerson \cite{Fu1} proved the
following min-max formula.

\THEOREM [Bock, Fulkerson] \label{fenyo-min-max} Let $c$ be a
non-negative cost function on the arc-set of a root-connected digraph
$D=(V,A)$.  \ The minimum cost of an $r_0$-arborescence of $D$ is
equal to \eq \hbox{ $\max \{\sum [ y(Z): \ Z\subseteq V-r_0] : \ y \ $is $
c$-feasible$\}$.  }\ \label{dualarbo} \eeq There is an optimal dual
solution $y$ for which $\{Z: \ y(Z)>0\}$ is laminar.  If $c$ is
integer-valued, the optimal $y$ can also be chosen integer-valued.
\FBOX \eT

Note that Fulkerson \cite{Fu1} developed a simple greedy algorithm for
computing the optimal dual vector $y$ occurring in the theorem.  The
theorem immediately implies the following optimality criteria.

\Corollary \label{arb-opt-krit} Let $y\sp *$ be a $c$-feasible
function on the family of non-empty subsets of $V-r_0$ and let $F_0$
be an $r_0$-arborescence of $D$ for which the following optimality
criteria hold:  $$ \hbox{ $F_0$ \ consists of tight arcs,}\ $$ $$
\hbox{ $y\sp *(Z)>0$ implies $\varrho _{F_0}(Z)=1$.}\ $$ \noindent
Then $F_0$ is a $c$-minimal $r_0$-arborescence of $D$ for which
$\widetilde c(F_0) = \sum [y\sp *(Z): \ Z\subseteq V-r_0]$.  \FBOX \eCo

In the following extension of Theorem \ref{fenyo-min-max}, a subset
$L\subseteq A$ of arcs is said to {\bf cover} the set-system $\cal F$
if $\varrho _L(Z)\geq 1$ holds for every $Z\in {\cal F}$.  The
set-system $\cal F$ is {\bf intersecting} if $X\cap Y, X\cup Y \in
{\cal F}$ holds whenever $X,Y\in {\cal F}$ and $X\cap Y \ne
\emptyset $.

\THEOREM [Frank \cite{FrankJ2}] \label{kernel-minmax} Let $\cal F$ be
an intersecting set-system on ground-set $V$ and $c:A\rightarrow {\bf
R} _+$ a cost function.  Then \eq \hbox{ $\min \{\widetilde c(L): \
L\subseteq A, \ L $ covers ${\cal F}\} = \max \{ \widetilde
y({\cal F}):$ \ $y:{\cal F}\rightarrow {\bf R} _+$ is
$c$-feasible$\}$.}\ \label{(kernel-minmax)} \eeq In addition, the
optimal $y$ can be chosen integer-valued when $c$ is integer-valued.
\FBOX \eT

We shall rely on the following fundamental result of Edmonds
\cite{Edmonds73}.

\THEOREM [Edmonds] \label{Edmonds} A digraph $D=(V,A)$ with root-node
$r_0$ has $k$ disjoint arborescences if and only if $D$ is rooted
$k$-arc-connected.\FBOX \eT

Note that Lov\'asz \cite{Lovasz76a} provided a particularly elegant
(algorithmic) proof for Theorem \ref{Edmonds}.

\section{The inverse cheapest arborescence problem }

Let $F_0$ be an $r_0$-arborescence of $D$ and let $A_0:=A-F_0$. We
shall refer to $F_0$ as the {\bf input} or {\bf target} arborescence.
We say that a cost function $w:A\rightarrow {\bf R} $ is an {\bf
$F_0$-minimizer} if $F_0$ is a $w$-minimal arborescence.  In the
inverse arborescence problem, we want to modify an input cost function
$w_0:A\rightarrow {\bf R} _+$ in such a way that the obtained cost
function $w$ is an $F_0$-minimizer, while the deviation of $w$ from
$w_0$ is minimized.  There are several ways to measure the
deviation:  here we focus on the $l_1$-norm.  Let $\vert z\vert :=
\vert z_1\vert +\cdots +\vert z_k\vert $ denote the $l_1$-norm of a
vector $z=(z_1,\dots ,z_k)$.

For an arborescence $F$, let $$\nabla(F):= \widetilde w_0(F_0) -
\widetilde w_0(F).$$ If $F\sp *$ is a minimum $w_0$-cost arborescence
of $D$, then $\nabla(F\sp *)$ is clearly a lower bound for the minimum
deviation we are looking for.  The following example shows that this
natural lower bound cannot always be achieved.

Let $V := \{r_0, u, v\}$ and $A := \{r_0u, r_0v, uv, vu\}$.  Let the
input cost function $w_0$ be $w_0(r_0u) = 1, \ w_0(r_0v) = 1, \
w_0(uv) = 0, \ w_0(vu) = 0$.  Let the target $r_0$-arborescence be
$F_0 := \{r_0u, r_0v\}$.  Now both $F_1 := \{r_0u, uv\}$ and $F_2 :=
\{r_0v, vu\}$ are cheapest $r_0$-arborescences with $w_0$-cost 1,
while $\widetilde w_0(F_0) = 2$.  It follows that the deviation of an
optimal $w$ from $w_0$ is at least $\widetilde w_0(F_0) - \widetilde
w_0(F_1) = 1$.  However, it is easy to check that, with a total
deviation of 1, it is not possible to obtain a cost function $w$ for
which $F_0$ is a $w$-minimal arborescence. 
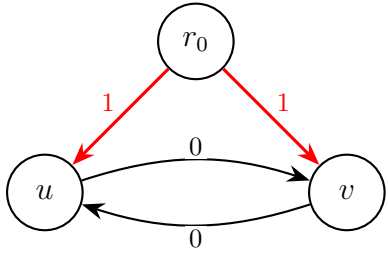
\begin{figure}[h]
\centering

\begin{tikzpicture}[
    vertex/.style={
        circle,
        draw,
        thick,
        minimum size=10mm,
        inner sep=0pt,
        font=\large
    },
    edge/.style={
        -{Stealth[length=3mm]},
        thick
    },
    target edge/.style={
        edge,
        red,
        very thick
    },
    edge label/.style={
        midway,
        fill=white,
        inner sep=1.5pt,
        font=\small
    }
]

\node[vertex] (r0) at (0,2) {$r_0$};
\node[vertex] (u)  at (-2,0) {$u$};
\node[vertex] (v)  at (2,0) {$v$};

\draw[target edge] (r0) -- node[edge label, above left] {$1$} (u);
\draw[target edge] (r0) -- node[edge label, above right] {$1$} (v);

\draw[edge] 
    (u) to[bend left=18] 
    node[edge label, above] {$0$} 
    (v);

\draw[edge] 
    (v) to[bend left=18] 
    node[edge label, below] {$0$} 
    (u);

\end{tikzpicture}
\caption{The natural lower bound cannot be achieved.}
\label{figure}
\end{figure}

\Remark {\em It may be useful to observe that in the analogous inverse
optimization problem concerning perfect matchings of a bipartite graph
$G$ the situation is simpler in the following sense.  If $F_0$ denotes
a target perfect matching of $G$ and $F\sp *$ is a minimum $w_0$-cost
perfect matching of $G$, then the minimum total change of $w_0$ (in
$l_1$-norm) to make $F_0$ a cheapest perfect matching is $\widetilde
w_0(F_0)- \widetilde w_0(F\sp *)$.  This can be proved rather easily
with the help of Egerv\'ary's foundational min-max theorem (see, for
example, Section 3.2.1 of the thesis of G. Hajdu \cite{Hajdu2020}).
$\bullet $ }\eRe

As far as the algorithmic aspect is concerned, Hu and Liu
\cite{Hu-Liu} developed a strongly polynomial algorithm for the
solution of the inverse arborescence problem.  They showed that there
always exists an optimal cost function $w\sp *$ which is non-negative
and {\bf $A_0$-fixed} in the sense that ${w\sp *}_{\vert A_0} =
{w_0}_{\vert A_0}$.  (These properties are not self-evident in the
light of the related inverse cheapest path problem, where they do not
hold.)  In addition, $w\sp *$ may be chosen integer-valued when $w_0$
is integer-valued.  Frank and Hajdu \cite{Frank+Hajdu} proved the following
min-max formula for the minimum deviation.


\THEOREM [\cite{Frank+Hajdu}] \label{dev-min-max} Let $w_0$ be a
non-negative cost function on the arc-set of digraph $D=(V,A)$ and let
$F_0$ be an $r_0$-arborescence of $D$.  Let ${\cal F}_0:= \{Z\subseteq
V-r_0:  \ \varrho _{F_0}(Z)=1\}$.  Then \eq \min \{\vert w-w_0\vert :
\ w \ \hbox{is an $F_0$-minimizer cost function$\} \ = $ }\
\label{(dev-min)} \eeq \eq \max \{\widetilde w_0(F_0)-\widetilde
w_0(L): \ L\subseteq A, \ L \ \hbox{covers ${\cal F}_0\}.$ }\
\label{(dev-max)} \eeq \noindent Moreover, there exists a
non-negative, $A_0$-fixed minimizer to \emeref{(dev-min)} which is
integer-valued when $w_0$ is integer-valued.  \FBOX \eT

The proof of the non-trivial direction $\max\geq \min$ relied on
Theorem \ref{kernel-minmax}.  The proof of this latter result in
\cite{FrankJ2} consists of a simple two-phase greedy algorithm, which
computes in its first phase an optimal $y$ occurring in Formula
\eref{(kernel-minmax)}, while an optimal $L$ is computed in the second
phase.  This approach, when applied in the special case ${\cal
F}:={\cal F}_0$, provided a direct algorithmic proof for the
non-trivial direction $\max\geq \min$ of Theorem \ref{dev-min-max}.

However, Theorem \ref{dev-min-max} has the (aesthetic) drawback that
the \lq trivial\rq \ direction $\max\leq \min$ is far from trivial, at
least in the sense that it relies on Theorem \ref{Edmonds} of Edmonds
(or, in an alternative approach, on the polyhedral description of the
polytope of $r_0$-arborescences which is implied by Theorem
\ref{fenyo-min-max}).

Our present goal is to formulate and prove a simplified min-max
theorem in which the inequality $\max \leq \min$ can indeed be seen
immediately.  In this new approach, it is the non-trivial direction
$\max \geq \min$ whose (algorithmic) proof shall rely on Edmonds'
theorem.


\subsection{A simplified min-max formula}

We say that a pair $\{F_1,F_2\}$ of arborescences {\bf covers $F_0$}
or that the pair is {\bf $F_0$-covering} if $F_0\subseteq F_1\cup
F_2$.  A simple observation for such a pair is that $F_1$ and $F_2$
are disjoint outside $F_0$, that is, $F_1\cap F_2\subseteq F_0$.
Indeed, if, indirectly, there is an arc $uv\in A-F_0$ belonging to
both $F_1$ and $F_2$, then the unique arc of $F_0$ entering $v$
belongs to at least one of $F_1$ and $F_2$, say $F_1$, implying that
the arborescence $F_1$ would have two arcs entering $v$.

\THEOREM Let $D=(V,A)$ be a digraph with a root-node $r_0\in V$.  Let
$F_0$ be an $r_0$-arborescence of $D$, and $w_0$ a non-negative cost
function on $A$.  Then

\eq \hbox{ $ \min \{ \vert w-w_0\vert :  \ w \ $is an $F_0$-minimizer cost function$\} \ =$ }\ \label{(min1)} \eeq \eq
\hbox{ $ \max \{ \nabla(F_1) + \nabla(F_2):  \ \{F_1,F_2\} \ \hbox{is
an}\ F_0$-covering pair of arborescences$\}$.  }\ \ \label{(max1)}
\eeq

\noindent There is a minimizer $w$ to \emeref{(min1)} which is
non-negative and $A_0$-fixed.  In addition, such a $w$ can be chosen
integer-valued when $w_0$ is integer-valued.  \eT

\Proof For a cost function $w$ minimizing \eref{(min1)}, one has
$w(a)\leq w_0(a)$ when $a\in F_0$, and $w(a)\geq w_0(a)$ when $a\in
A_0$.  For an arbitrary arborescence $F$, the value
$\nabla(F)=\widetilde w_0(F_0) - \widetilde w_0(F)$ provides a lower
bound for the total change of $w_0$ on $F_0\oplus F$. By applying
this observation separately to the two members of an $F_0$-covering
pair $\{F_1,F_2\}$ of arborescences, the disjointness of the sets
$F_0\oplus F_1$ and $F_0\oplus F_2$ imply that \eq
\vert w-w_0\vert \geq [\widetilde w_0(F_0) - \widetilde w_0(F_1)] +
[\widetilde w_0(F_0) - \widetilde w_0(F_2)] = \nabla(F_1) +
\nabla(F_2), \label{(triv)} \eeq from which $\max\leq \min$ follows.

To prove the reverse direction, we construct an $F_0$-minimizer cost
function $w\sp *$ which is non-negative and $A_0$-fixed, along with an
$F_0$-covering pair $\{F_1\sp *,F_2\sp *\}$ of arborescences for which
\eref{(triv)} holds with equality.

Apply Theorem \ref{kernel-minmax} to the intersecting set-system
${\cal F}_0:=\{X\subseteq V-r_0:  \ \varrho _{F_0}(X)=1\}$ in place of
$\cal F$ and to the cost function $c:=w_0$.  Let $L\sp * \subseteq A$ be a
minimum $c$-cost arc-set covering ${\cal F}_0$ which is minimal with
respect to inclusion, and let $y\sp * :  {\cal F}_0 \rightarrow {\bf
Z} _+$ be an optimal dual solution.

\Claim \label{1-befoku} $\varrho _{L\sp *}(v)=1$ holds for every node
$v\in V-r_0$.  \eCl

\Proof Since $\varrho _{F_0}(v)=1$, we have $\{v\}\in {\cal F}_0$ from
which $\varrho _{L\sp *}(v)\geq 1$.  Suppose indirectly that there are
two arcs $a_1$ and $a_2$ in $L\sp *$ entering $v$.  The minimality of
$L\sp *$ implies that ${\cal F}_0$ has a member $Z_i$ $(i=1,2)$ for
which $a_i$ is the single element of $L\sp *$ entering $Z_i$.  Since
$v\in Z_1\cap Z_2$, both the intersection $Z_1\cap Z_2$ and the union
$Z_1\cup Z_2$ are in ${\cal F}_0$, and hence the submodularity of
$\varrho _{L\sp *}$ implies that $$1+1= \varrho _{L\sp *}(Z_1) +
\varrho _{L\sp *}(Z_2) \geq \varrho _{L\sp *}(Z_1\cap Z_2) + \varrho
_{L\sp *}(Z_1\cup Z_2) \geq 1+1,$$ from which $\varrho _{L\sp
*}(Z_1\cap Z_2)=1$, contradicting the property that both $a_1$ and
$a_2$ enter $Z_1\cap Z_2$.  \FBOX

\medskip

Let $F_0\sp *$ denote the arborescence whose arcs are parallel copies
of the arcs of $F_0$ (and hence $F_0\sp *$ is disjoint from $A$).  Now
digraph $D':=(V,L\sp * \cup F_0\sp *)$ is rooted 2-arc-connected since
if $Z\subseteq V-r_0$ is a set entered by a single element of $F\sp
*_0$, then $Z\in {\cal F}_0$, and $L\sp *$ contains an arc entering
$Z$.  By Claim \ref{1-befoku}, the in-degree of each node $v\in
V-r_0$ in $D'$ is $2$, and hence Theorem \ref{Edmonds} of Edmonds
implies that the arc-set of $D'$ can be partitioned into two spanning
arborescences of root $r_0$.

Let $F_1\sp *$ and $F_2\sp *$ denote the arborescences of $D$
corresponding to these two arborescences of $D'$.  It follows from
this construction that the pair $\{F_1\sp *, F_2\sp *\}$ of
arborescences is $F_0$-covering, for which $L\sp *=F_1\sp * \cup
F_2\sp *$.  Let \eq w\sp *(a):= \begin{cases} c(a) \ \ \ (=w_0(a)) & \
\ \hbox{if}\ \ \ a\in A_0 \cr \sum [y\sp *(Z) : \ a \ \hbox{enters} \
Z] & \ \ \hbox{if}\ \ \ a\in F_0.  \end{cases} \label{(w*def)} \eeq

\noindent This $w\sp *$ is non-negative and $A_0$-fixed.  Furthermore,
$y\sp *$ is feasible with respect to $w\sp *$ as well, and every arc
in $F_0$ is $w\sp *$-tight.  By applying Corollary \ref{arb-opt-krit}
to $c:=w\sp *$, we obtain that $F_0$ is a $w\sp *$-minimal
arborescence.

It follows from \eref{(kernel-minmax)} that $w_0(a) = \sum [y\sp
*(Z):  \ a$  enters $Z]$ holds for every arc $a\in L\sp *$.
Therefore, $w_0(a) = w\sp *(a)$ holds for every element $a\in L\sp *
\cap F_0 \ (=F\sp *_1\cap F\sp *_2)$.  Since $w_0(a) = w\sp *(a)$ also
holds for the elements of $A_0$, we get $$\widetilde w_0(L\sp
*) = \widetilde w\sp *(L\sp *) = \widetilde y\sp *({\cal F}_0) =
\widetilde w\sp *(F_0),$$ from which $$ \vert w\sp *-w_0\vert =
\widetilde w_0(F_0) - \widetilde w\sp *(F_0) = \widetilde w_0(F_0) -
\widetilde w_0(L\sp *) = \nabla (F_1\sp *) + \nabla(F_2\sp *),$$
showing that the cost function $w\sp *$ and the $F_0$-covering pair
$\{F_1\sp *, F_2\sp *\}$ of arborescences indeed satisfy
\eref{(triv)} with equality.  \BB

\medskip

As mentioned above, the proof of Theorem \ref{kernel-minmax} occurring
in \cite{FrankJ2} is a two-phase greedy algorithm which computes an
optimal dual solution $y\sp *$ in the first phase and an optimal
primal solution $L\sp *$ covering ${\cal F}_0$ in the second phase.
Once $L\sp *$ is available, one can apply the algorithmic proof of the
disjoint arborescence theorem of Edmonds \cite{Edmonds73} (in the
special case of $k=2$) provided by Lov\'asz \cite{Lovasz76a},
obtaining in this way the optimal $F_0$-covering pair $\{F_1\sp
*,F_2\sp *\}$ of arborescences occurring in \eref{(max1)}.

It should be noted that the paper \cite{Frank+Hajdu} of Frank and
Hajdu described a specific, more efficient two-phase greedy algorithm
for the intersecting set-system ${\cal F}_0$ introduced in Theorem
\ref{dev-min-max}.


\medskip

\begin{thebibliography}{999}




\bibitem{Bock} F. Bock, {\it An algorithm to construct a minimum
directed spanning tree in a directed network}, in:  Developments of
Operations Research, Vol. 1. (Proceedings of the Third Annual Israel
Conference on Operations Research, Tel Aviv, 1969, B. Avi-Itzhak, ed.)
Gordon and Breach, New York, 1971, 29--44.


\bibitem{Edmonds73} J. Edmonds, {\it Edge-disjoint branchings,} in:
Combinatorial Algorithms (B.  Rustin, ed.), Acad.  Press, New York,
(1973) 91-96.




\bibitem{FrankJ2} A. Frank, {\it Kernel systems of directed graphs},
Acta Scientiarum Mathematicarum (Szeged) 41, 1-2 (1979) 63-76.


\bibitem{Frank+Hajdu} A. Frank and G. Hajdu, {\it A simple algorithm
and min-max formula for the inverse arborescence problem}, Discrete
Applied Mathematics, 295 (2021) 85--93.




\bibitem{Fu1} D.R.  Fulkerson, {\it Packing rooted directed cuts in a
weighted directed graph}, Math.  Programming 6 (1974) 1-13.


\bibitem{Hajdu2020} G. Hajdu, {\it Inverz kombinatorikus
optimaliz\'al\'asi probl\'em\'ak} \ (Inverse combinatorial
optimization problems), BSc thesis in Hungarian (2020), E\"otv\"os
Lor\'and University Budapest.

\bibitem{Hu-Liu} Z. Hu and Z. Liu, {\it A strongly polynomial
algorithm for the inverse shortest arborescence problem}, Discrete
Applied Mathematics, 82 (1998) 135--154.



\bibitem{Lovasz76a} L. Lov\'asz, {\it On two minimax theorems in graph
theory,} J. Combinatorial Theory (B), 21 (1976) 96-103.

\end{thebibliography}
\end{document}